\documentclass [12pt] {amsart}
\usepackage{palatino}
\usepackage{latexsym,amssymb, amsmath}
\usepackage{bm}
\usepackage{setspace}

\usepackage{simplemargins}
\setleftmargin{2.75cm}
\setrightmargin{2.75cm}
\settopmargin{2.0in}
\setbottommargin{2.0in}
\setallmargins{2.75cm}

\def\pf{\noindent\emph{Proof: }}
\def\stop{\hfill$\Box$}

\usepackage{amsmath, amsfonts, amssymb, amsthm}
\newtheorem{thm}{Theorem}

\newtheorem{lemma}[thm]{Lemma}

\newtheorem{defn}[thm]{Definition}

\numberwithin{thm}{section}

\DeclareMathOperator{\Hol}{Hol}

\begin{document}

\title [On Local Solutions of Second Order Quasi-linear Elliptic Systems] {On Local Solutions of Second Order Quasi-linear Elliptic Systems with Arbitrary 1-Jet at a Point}

\author{Yifei Pan}
\address{Department of Mathematical Sciences\\
 Purdue University Fort Wayne\\
 Fort Wayne, Indiana 46805}
\email{pan@pfw.edu}

\author{ Yu Yan}
\address {Department of Mathematics and Computer Science\\ Biola University\\La Mirada, California 90639}
\email {yu.yan@biola.edu}


\begin{abstract}
We prove a general result on the existence of local solutions of any second order quasi-linear elliptic system with arbitrary 1-jet at a point.
\end{abstract}

\maketitle

\section {Introduction}

\vspace{.1in}
A fundamental question in the study of PDEs is the existence of solutions.  It is well understood for linear PDEs; for nonlinear equations there have been numerous studies of specific equations. In this paper we prove the existence of local solutions of general second order quasi-linear elliptic systems.  This result can be applied to many well known geometric equations and be used to define a new Kobayashi metric on Riemannian manifolds.

\vspace{.1in}



\vspace{.1in}

Let $\displaystyle \bm{u}(x)=\big ( u^1(x) , u^2(x), ..., u^m(x)\big ) $ be a $C^2$ vector function from an open set in $\mathbb{R}^n (n \geq 2)$ to $\mathbb{R}^m (m \geq 1)$. We use $L$ to denote an elliptic quasi-linear operator on $\bm{u}$, and write 
$$ L\bm{u}= \left ( Lu^1,..., Lu^m  \right ),$$ where
\begin{equation*}
	\label{defn:L}
	Lu^k= \sum _{i,j=1}^na^{ij}\big ( x, \bm{u}(x), D\bm{u}(x)\big ) D_{ij} u^k, \hspace{.2in} 1 \leq k \leq m,
\end{equation*}
\noindent
$a^{ij} \in C_{loc}^{1,\alpha}\left (\mathbb{R}^n \times\mathbb{R}^m \times \mathbb{R}^{mn} \right )$ ($0 < \alpha <1 $), and there is a $\lambda >0$ such that

$$a^{ij}(x,p,q)\xi_i\xi_j\geq \lambda|\xi|^2 $$
for all $\xi\in\mathbb{R}^n$.

\vspace{.1in}
\noindent
We will prove that any quasi-linear system defined by such an operator $L$ is always locally solvable, given the value and the first derivatives of the solution at any one point in $\mathbb{R}^n$.  Without loss of generality we can assume the domain of $\bm{u}$ is a ball centered at the origin.

\vspace{.1in}

\noindent
Denote a ball centered at the origin with radius $R$ in $\mathbb{R}^n$ as
$$ B_R = \{ x \in \mathbb{R}^n \big | |x| \leq R \} $$ and denote a ball centered at the origin in $\mathbb{R}^m$ as $$ B'_{R'}=\{ y \in \mathbb{R}^m \big | |y| \leq R' \}.  $$  The following is our main result.

\vspace{.1in}

\begin{thm}
	\label{thm: main}
	Let $\bm{\phi} (x, p, q)= \big (\phi ^1(x, p, q), ..., \phi ^m(x, p, q)\big ):B_R \times B'_{R'} \times \mathbb{R}^{mn} \to \mathbb{R}^{m}  $ be of $C_{loc}^{1,\alpha}$ $ (0 < \alpha < 1)$.
	For any given $\bm{c_0} \in B'_{R'}$ and $ \bm{c_1}\in \mathbb{R}^{mn} $, the elliptic quasi-linear system
	\begin{equation}
		\label{eqn:main}
		\renewcommand{\arraystretch}{1.5}
		\left \{
		\begin{array}{r@{}l}
		L\bm{u}(x) &{} = \bm{\phi}(x, \bm{u}(x), D\bm{u}(x)) \\
	\bm{u}(0) &{} =  \bm{c_0} \\
		D\bm{u}(0) &{} =  \bm{c_1}
	\end{array}
\right.
	\end{equation}

\noindent
	has $C^{2,\alpha}$ solutions $\bm{u}(x)$ from $B_R$ to $B'_{R'}$ when $R$ is sufficiently small.
	
\end{thm}

\noindent
As will be shown in the proof, there are in fact infinitely many such solutions.  An interesting feature of our theorem is that the only assumptions about the coefficient functions are ellipticity and sufficient regularity, therefore it is as general as we can hope for.  It is intriguing to know if such a general existence theorem holds for fully nonlinear elliptic systems. The equation $e^{\Delta u} =0 $, although being non-elliptic, seems to suggest that general local existence may not be true in the fully nonlinear case.

\vspace{.2in}

Theorem \ref{thm: main} has a wide range of applications in differential geometry, since many of the well known geometric equations are semilinear or quasilnear.

\vspace{.15in}

\noindent
\textbf{Example 1.} 
	(\textbf{The Minimal Surface Equation})  For a domain $\Omega \subset \mathbb{R}^n$ and a smooth function $u$ on $\Omega$, the area of its graph $(x, u(x)), $ where $x \in \Omega$, is $\displaystyle \int _{\Omega} \sqrt{1+|Du|^2} \,\,dx.  $  Its Euler-Lagrange equation is the minimal surface equation 
	$$ \sum_{i,j=1}^n   \left (\delta_{ij}-  \frac{D_iu D_ju}{1+|Du|^2}  \right )D_iD_ju=0, $$ which is an elliptic quasilinear equation. The graph of $u$ is called a minimal surface, on which the mean curvature is identically 0.  Without loss of generality we can assume $0\in \Omega$. By Theorem \ref{thm: main}, given $u(0)$ and $Du(0)$, this equation always has local solutions near $0$.  Therefore, given any point  $P=(0,...,0, p) \in \mathbb{R}^{n+1}$ and any normal vector $\eta$, there are infinitely many local minimal surfaces through $P$ with $\eta$ as the normal vector. 
	
	\vspace{.1in}
	\noindent
	Similar existence also holds for the prescribed mean curvature equation $$ \sum_{i,j=1}^n   \left (\delta_{ij}-  \frac{D_iu D_ju}{1+|Du|^2}  \right )D_iD_ju=H(x,u), $$ where $H(x,z)$ is a given function on $\Omega \times \mathbb{R}$.

\vspace{.2in}





\noindent
\textbf{Example 2.}
	 (\textbf{The Harmonic Map System}) Let $(M_1^n, g)$ and $(M_2^m,h)$ be two Riemannian manifolds of dimensions $n$ and $m$ with metrics $g$ and $h$, respectively.  A map $u:M_1^n \to M_2^m$ is called a harmonic map if it is a critical point of the energy functional. Let $x=(x^1,..., x^n)$ be a local coordinate system on $M_1$. In a local coordinate system on $M_2$, we may write $u(x)=(u^1(x),..., u^m(x))$, and it satisfies the quasi-linear equations
	 
	 	$$ \Delta u^{\alpha} +  \sum_{i,j=1}^n \sum _{\beta,\gamma=1}^m g^{ij}\Gamma_{\beta \gamma}^{\alpha}(u) \frac{\partial u^{\beta}}{\partial x^i} \frac{\partial u^{\gamma}}{\partial x^j} = 0, \hspace{.1in} \alpha =1,... , m,$$ 
	 	where $\Gamma_{\beta \gamma}^{\alpha}  $ are the Christoffel symbols on $M_2$.
	
	\vspace{.1in}
	\noindent
	By Theorem \ref{thm: main}, given any $p_1 \in M_1$, $p_2 \in M_2$, and a subspace $\Gamma$ of the tangent space $T_{p_2}M_2$, there exist infinitely many local harmonic maps $u: M_1 \to M_2$ which map $p_1$ to $p_2$ and $T_{p_1}M_1$ to $\Gamma$. The theorem also implies the existence of any local objects to be defined by the Laplace-Beltrami operator.

\vspace{.1in}

A novel application of this existence of harmonic maps is a definition of Kobayashi metric on Riemmannian manifolds.  The Kobayashi metric was first introduced as a pseudometric on complex manifolds by Kobayashi \cite{Kobayashi}.  Let $D=\{ \xi \in \mathbf{C} : |\xi| <1  \}$ be the unit disc in the complex plane, and let $M$ be a complex manifold. We denote the set of all holomorphic functions from $D$ to $M$ as $\Hol(D,M)$. For any $p \in M$ and $v \in T_p M$, the infinitesimal Kobayashi metric $K_M$ of $v$ at $p$ is defined by
\begin{eqnarray*}
	 K_M (p, v)
	 & = & \inf \left \{ \alpha : \,\, \alpha >0 \,\,  \text{and} \,\, \exists \,\, f \in \Hol(D,M) \,\,  \text{with} \,\, f(0)=p, f'(0)= \frac{v}{\alpha}  \right  \}  \\
	& = &	\inf \left \{ \frac{|v|}{|f'(0)|} \,\, : \,\, \exists \,\, f \in \Hol(D,M) \,\,  \text{with} \,\, f(0)=p  \right  \}  
\end{eqnarray*}

\vspace{.1in}
\noindent
It can be proved that $K_M$ is upper semicontinuous on $TM$, and therefore it is a Finsler metric on $M$.
The Kobayashi metric is a well studied, important tool in several complex variables and complex geometry, so a natural question is whether a similar concept can be introduced on Riemannian manifolds. Theorem \ref{thm: main} paves the way for answering this question by establishing a well-defined, real version of Kobayashi metric as follows.

\vspace{.1in}

\noindent
 Again we use $D$ to denote the unit disk in $\mathbb{R}^2$, and let $D_R$ be a disk of radius $R$ in $\mathbb{R}^2$ centered at the origin with coordinates $(x,y$).  Let $(M,g)$ be a Riemannian manifold with metric $g$. We consider harmonic maps  $u: D_R \to M$ that are conformal  at $0$, so  $$g \left (\frac{\partial u}{\partial x}(0),\frac{\partial u}{\partial x}(0) \right ) = g \left (\frac{\partial u}{\partial y}(0),\frac{\partial u}{\partial y}(0) \right ) \hspace{.1in} \text{and} \hspace{.1in} g \left (\frac{\partial u}{\partial x}(0),\frac{\partial u}{\partial y}(0) \right )=0.$$


\begin{defn}
	\label{defn:kobayashi}
	Let $(M,g)$ ba a Riemannian manifold, and let $X \in T_p M$ be a non-zero tangent vector at a point $p \in M$.  The Kobayashi metric of $X$ is defined by		
	\begin{eqnarray}
		\label{eqn:kobayashi}	
	& &	K_M(p,X) \\
		& =&   \inf \left \{ \alpha: \alpha >0, \,\, \exists \text{ harmonic map } u: D \to M, \,\,\text{conformal at }  0,\,\, u(0)=p, \,\, \frac{\partial u}{\partial x}(0)= \frac{X}{\alpha} \right \} \nonumber \\
	& =&   \inf \left \{ \frac{1}{R}: \,\, \exists \text{ harmonic map } u: D_R \to M, \,\, \text{conformal at }  0,\,\, u(0)=p, \,\, \frac{\partial u}{\partial x}(0)=X \right \} \nonumber 
	\end{eqnarray}
\noindent
If $X=0$, we define $K_M(p,0)=0$.

\end{defn}
	
\vspace{.1in}

\noindent
For any non-zero vector $X \in T_p M$, we can find another vector $Y$ satisfying $g(X,X)=g(Y,Y)$ and $ g(X,Y) =0$.  Then by Theorem \ref{thm: main}, there is a harmonic map $u$ on a local disk $D_R, R>0$, such that $u(0)=p$, $\displaystyle \frac{\partial u}{\partial x}(0)=X$, and $\displaystyle \frac{\partial u}{\partial y}(0)=Y$.  This function $u$ is conformal at $0$ by the choice of $X$ and $Y$. Therefore in (\ref{eqn:kobayashi}) the infimum of $\displaystyle \frac{1}{R}$ is always a finite number, and consequently the Kobayashi metric is well defined on Riemannian manifolds:

\vspace{.1in}

\begin{thm}
	Let $(M,g)$ be any Riemannian manifold.  The Kobayashi metric in (\ref{eqn:kobayashi}) is well defined.  Namely, $0 \leq K_M(p,X) < \infty $ for all $X \in T_p M $.
\end{thm}

\vspace{.1in}
\noindent
We would like to point out that the condition of conformality at 0 excludes geodesics from the admissible harmonic maps.  For any geodesic $\gamma(t)$ on $M$ such that $\gamma(0)=p$ and $ \gamma '(0)=X$, the function $\displaystyle u(x,y)=\gamma (x+y)$ defines a harmonic map from a neighborhood of the origin to $M$, with $u(0)=p$ and $\displaystyle \frac{\partial u}{\partial x}(0)=X$. However, since $\displaystyle \frac{\partial u}{\partial y}(0)=X$ too, $u$ is not conformal at $0$.

  
\vspace{.1in}

There are a lot of questions we can ask about $K_M$.  The first one is whether it is a Finsler metric, which needs to be upper-semicontinuous and positive definite.  We hope to explore this in future work.  


\vspace{.1in}




\vspace{.2in}
The rest of the paper is organized as follows. The main strategy for proving Theorem \ref{thm: main} is similar to that in \cite{Pan_Zhang}.  The solution is found by applying the Fixed Point Theorem to an appropriately defined Banach space, but the quasilinear term introduces additional subtlety that needs to be handled carefully.  In Section \ref{sec: Banach} we define a Banach space of functions with vanishing order from which we will seek possible solutions. Then, we study the Newtonian potential as an operator on this Banach space in Section \ref{sec:potential}. Next, in Section \ref{sec:Poisson} we show that proving Theorem \ref{thm: main} is equivalent to proving Lemma \ref{thm: Poisson}, which is the local existence of solutions of a Poisson type system, and we define a map for the Poisson system.  Finally, in Section \ref{sec:estimate_f} and Section \ref{sec:estimate_f-g} we prove this map is a contraction if some parameters are appropriately chosen, then the Fixed Point Theorem can be invoked to find a solution, proving Lemma \ref{thm: Poisson}.  Thus, Theorem \ref{thm: main} is also true.

\vspace{.2in}

\section{Defining the Function Spaces and Some Preliminary Estimates}
\label{sec: Banach}

\vspace{.1in}

First we define some notations we will be using throughout the proof.  

\vspace{.1in}

\noindent
Let $f: B_R \to \mathbb{R}$ be a function defined on $B_R$.

\begin{itemize}
	
	\item 
	The H\" older seminorm of $f$ is $$\displaystyle H_{\alpha}[f]=\sup _{x, x' \in B_R}\bigg\{\frac{|f(x)-f(x')|}{|x-x'|^{\alpha}} \bigg|  \bigg\}.$$
	
	\vspace{.1in}
	
	\item $C^{0,\alpha}(B_R)$ is the set of all functions on $B_R$ such that $H_{\alpha}[f]$ is finite.
	
	\vspace{.1in}
	
	\item The weighted H\" older norm of $f$ is 
	\begin{equation}
		\label{def:weighted_norm}
	\displaystyle \|f\|_{\alpha}=\sup_{x\in B_R}|f(x)|+(2R)^\alpha H_\alpha[f].
	\end{equation} 
	It is well known (see \cite{GT}) that $\| \cdot \|_{\alpha}$ is a norm on $C^{0,\alpha}(B_R)$ with which  $C^{0,\alpha}(B_R)$ is a Banach algebra, i.e.	
	\begin{equation}
		\label{eqn:Banach_algebra}
		\| fg \|_{\alpha} \leq \|f\|_{\alpha} \| g \|_{\alpha}
	\end{equation}
	for any $f,g \in C^{\alpha}(B_R)$.
	
	\vspace{.1in}
	
	\item The H\" older space $C^{2,\alpha}(B_R)$ consists of all functions $f$ on $B_R$ whose second order partial derivatives exist and belong to $C^{0,\alpha}(B_R)$.
	
	\vspace{.1in}
	
	\item We use $C_0^{2,\alpha}(B_R)$ to denote the set of all functions in $C^{2,\alpha}(B_R)$ whose first order derivatives all vanish at the origin:
	$$C_0^{2,\alpha}(B_R)=  \{f\in  C^{2,\alpha}(B_R)\big | \partial_i f(0)=0, \,\,i=1,..., n\}.$$
	
	\vspace{.1in}
	\item For any integer $l \geq 0$, define 
	\allowdisplaybreaks
	\begin{eqnarray}
		\label{eqn:(k)-norm}
		\|f\|^{(l,\alpha)} & = &  \max_{|\beta|=l}  \big \{\|\partial^\beta f\|_{\alpha} \big \} \nonumber \\
		 & = &  \max_{|\beta|=l}  \Big \{\sup_{x\in B_R}|\partial ^{\beta}f(x)|+(2R)^\alpha H_\alpha[\partial ^{\beta}f]\Big \}, 
\end{eqnarray}
	where we have used the notation $\beta=(\beta_1,...,\beta_n)$, $|\beta|=\beta_1+...+\beta_n$, and	$$\partial^\beta=\partial_1^{\beta_1}\cdot\cdot\cdot\partial_n^{\beta_n}.$$ Note that $\|f\|^{(0,\alpha)} = \|f\|_{\alpha}. $

\vspace{.1in}

\noindent
We can also extend the definition $\|\cdot \|^{(l, \alpha)}$ to vector functions  $\bm{f}=(f^1,..., f^m): B_R \to \mathbb{R}^m$ and define
\begin{equation}
	\label{defn:vector_norm}
	\|\bm{f}\|^{(l,\alpha)} = \max _{j=1,..., m} \|f^j \|^{(l,\alpha)}.
\end{equation}

\end{itemize}
 
\vspace{.1in}
In our proof we will only use (\ref{eqn:(k)-norm}) and (\ref{defn:vector_norm}) for $0 \leq l \leq 2$.	A key fact to be used in our proof is that $\|\cdot\|^{(2, \alpha)}$ is a norm on $C_0^{2,\alpha}(B_R)$, and $C_0^{2,\alpha}(B_R)$ becomes a Banach space under this norm.  We will establish this fact in Lemma \ref{lem:C_0_Banach}, which is based on the estimates in the following Lemmas \ref{lem:remainder} to \ref{lem: compare_norms}.

\vspace{.1in}

\begin{lemma} 
	\label{lem:remainder}	
	If $f\in C^{2,\alpha}(B_R)$, then for any $x, x'\in B_R$,
	$$\Big |f(x')-\sum_{l=0}^2  \frac{1}{l!} \sum_{|\beta|=l} \partial^\beta f(x)(x'-x)^\beta \Big |\leq \frac{1}{2} \left ( \sum_{|\beta|=2} H_\alpha[\partial^\beta f]\right )|x'-x|^{2+\alpha}.$$
\end{lemma}

\pf  Define $$\varphi (t) = f ( x+th ), \hspace{.2in} \text{where}  \hspace{.2in} h=x'-x.$$  Then

\allowdisplaybreaks
\begin{eqnarray*}
	\int _0 ^1  \int _0 ^ {t_1} \varphi ^{''} (t) dt dt_1 	& = & \varphi (1) - \varphi (0) - \varphi ' (0) \\
	& = & f(x') - f(x) - \sum _{i=1}^n \partial _i f(x)h_i  ,
\end{eqnarray*}



\noindent
 After subtracting the second derivatives term, we have	 
 
 \allowdisplaybreaks
 \begin{eqnarray*}
 & & \Big | 	f(x') - f(x) - \sum _{i=1}^n \partial _i f(x)h_i -\frac{1}{2} \sum _{i,j=1}^n \partial _{i} \partial_j f(x)h_ih_j   \Big |   \\
 	& = &   \Bigg |  \int _0 ^1 \int _0 ^ {t_1} \varphi ^{''} (t) dt dt_1  - \frac{1}{2} \sum _{i,j=1}^n \partial _{i} \partial_j f(x)h_ih_j  \Bigg |   \\
 	& = &   \Bigg |  \int _0 ^1  \int _0 ^ {t_1} \varphi ^{''} (t) dt dt_1   - 	\int _0 ^1  \int _0 ^ {t_1} \sum _{i,j=1}^n \partial _{i} \partial_j f(x)h_ih_j  dt dt_1  \Bigg |  \\
 	& = &  \Bigg |  \int _0 ^1 \int _0 ^ {t_1} \left ( \sum _{i,j=1}^n \partial _{i} \partial_j f(x+th)h_ih_j -\sum _{i,j=1}^n \partial _{i} \partial_j f(x)h_ih_j  \right ) dt dt_1  \Bigg |  \\
 	& \leq & \int _0 ^1  \int _0 ^ {t_1}  \sum _{i,j=1}^n  H_{\alpha} \big [ \partial _{i} \partial_j f  \big ] |h|^{2+\alpha}     dt dt_1 \\
 	& = & \frac{1}{2}\sum _{i,j=1}^n  H_{\alpha} \big [ \partial _{i} \partial_j f  \big ] |h|^{2+\alpha},
 \end{eqnarray*}


\noindent
which can be written as 
	$$\Big |f(x')-\sum_{l=0}^2  \frac{1}{l!} \sum_{|\beta|=l} \partial^\beta f(x)(x'-x)^\beta \Big |\leq  \frac{1}{2} \left ( \sum_{|\beta|=2} H_\alpha[\partial^\beta f]\right )|x'-x|^{2+\alpha}.$$

\stop

\begin{lemma}
	\label{lem:f-bounded-by_f(k)}
	If $f\in C_0^{2,\alpha}(B_R)$, then
	$$\|f\|_{\alpha} < (3nR)^2\|f\|^{(2, \alpha)}.$$
\end{lemma}

\pf
	Let $f\in C_0^{2,\alpha}(B_R)$, by definition $f(0)=0$ and $\displaystyle \partial _i f(0)=0$ for all $i= 1,...,n $.  
	
\vspace{.1in}	
	\noindent
	For any $x \in B_R$, define	
 $$\varphi(t)=f(tx),$$ then $\varphi(0)=0$ and $\varphi '(0)=0$.  Therefore,
	
	 \allowdisplaybreaks
	\begin{eqnarray}
		\label{eqn:f_use_Banach_algebra}
		f(x)&=&\int_0^1 \int_0^{t_1} \varphi ''(t) dtdt_1 \nonumber \\
		&=&\int_0^1 \int_0^{t_1}\left ( \sum_{|\beta|=2}\partial^\beta f(t x)x^\beta\right )  dtdt_1 \nonumber \\
		&=&\sum_{|\beta|=2} \left ( \int_0^1 \int_0^{t_1}\partial^\beta f(tx) \,\, dtdt_1 \right ) x^\beta.
	\end{eqnarray}

\vspace{.1in}
\noindent
It follows easily from (\ref{def:weighted_norm}) that $\|x_i\|_{\alpha}=3R$ for all $ i=1,...n$.  Then for any $\beta = (\beta_1,..., \beta_n)$ with $|\beta|=2$, by (\ref{eqn:Banach_algebra}) we have 
$$	\|x^\beta\|_{\alpha}  =   \|x_1^{\beta _1} \cdots x_n^{\beta _n} \|_{\alpha} \\
	 \leq  (3R)^2. $$

\noindent
Now applying (\ref{eqn:Banach_algebra}) to (\ref{eqn:f_use_Banach_algebra}), we obtain	
	 \allowdisplaybreaks
	\begin{eqnarray*}
		\|f\|_{\alpha} &\leq & \sum_{|\beta|=2} \|\partial^\beta f\|_{\alpha}\|x^\beta\|_{\alpha} \\
		& \leq & (3R)^2 \sum_{|\beta|=2} \|\partial^\beta f\|_{\alpha}  \\
	&	\leq & (3nR)^2\|f\|^{(2,\alpha)}.
	\end{eqnarray*}

\stop

\begin{lemma}
	\label{lem: compare_norms}
	If $f\in C_0^{2,\alpha}(B_R)$, then for any $l < 2$,
	$$\| f \|^{(l,\alpha)}  \leq (3nR)^{2-l} \|f\|^{(2, \alpha)}.$$
\end{lemma}

\pf
	Let  $f\in C_0^{2,\alpha}(B_R)$. If $|\beta|=l < 2$, then $\partial^\beta f\in C_0^{2-l ,\alpha}(B_R)$. 
	
	\vspace{.1in}
	\noindent
	Similar to the proof of Lemma \ref{lem:f-bounded-by_f(k)}, we can show that
	\allowdisplaybreaks
\begin{eqnarray*}
	\|\partial^\beta f\|_{\alpha} & \leq & (3nR)^{2-l}\|\partial^\beta f\|^{(2-l, \alpha)} \\
	& \leq &  (3nR)^{2-l} \| f\|^{(2, \alpha)}.
\end{eqnarray*}

\noindent
Therefore $$\| f \|^{(l,\alpha)}  \leq (3nR)^{2-l} \|f\|^{(2, \alpha)}.$$
\stop

\vspace{.2in}

Next, we prove the main result of this section.
	\vspace{.1in}
	
\begin{lemma}
	\label{lem:C_0_Banach}
	The function space $C_0^{2,\alpha}(B_R)$ equipped with the  $\|\cdot\|^{(2, \alpha)}$ norm is a Banach space.
\end{lemma}
	
\pf  By definition $\| \cdot \|^{(2, \alpha)} $ is a semi-norm on $C^{2,\alpha}(B_R)$.  If $\| f \|^{(2, \alpha)}=0$ , then $\partial ^{\beta} f =0$ for all $|\beta| =2$ on $B_R$, which implies $f$ is a constant or a linear function.  If in addition $f \in C_0^{2,\alpha}(B_R)$, then $f$ and its first derivatives all vanish at $0$, so $f$ must be identically $0$.  Thus $\| \cdot \|^{(2, \alpha)} $ is a norm on $C_0^{2,\alpha}(B_R)$. 

\vspace{.1in}
\noindent
Since $C_0^{2,\alpha}(B_R)$ is a closed subspace of $C^{2,\alpha}(B_R)$ and $C^{2,\alpha}(B_R)$ is a Banach space with the $\| \cdot \|_{\alpha} $ norm, we know $C_0^{2,\alpha}(B_R)$ is a Banach space with the $\| \cdot \|_{\alpha} $ norm.  Then by Lemma \ref{lem: compare_norms}, $C_0^{2,\alpha}(B_R)$ is also complete under the $\| \cdot \|^{(2, \alpha)} $ norm.  Therefore, $C_0^{2,\alpha}(B_R)$ equipped with the  $\|\cdot\|^{(2, \alpha)}$ norm is a Banach space. 

\stop

\vspace{.2in}

\section{A H\" older Estimate for the Newtonian Potential}
\label{sec:potential}

\vspace{.1in}

Recall that the fundamental solution of the Laplace's equation is given by

$$ \Gamma(x-y)= 
\begin{cases}
	\displaystyle \frac{1}{n(n-2)\omega_n}|x-y|^{2-n} & \text{if} \,\,n \geq 3, \\
\displaystyle -	\frac{1}{2\pi}\ln|x-y| & \text{if} \,\, n=2. 
\end{cases} $$

\vspace{.1in}

\noindent
For an integrable function $f$ on $B_R$, the Newtonian potential of $f$ is  defined on $\mathbb{R}^n$ by
$$\mathcal{N}(f)(x)=\int_{B_R} \Gamma(x-y)f(y)dy;$$ it solves the Poisson's Equation $\Delta u=- f$.  In our proof we will consider $\mathcal{N}$ as an operator acting on a function space. The following result is well-known and a proof is given in \cite{Fraenkel}.

\begin{lemma} 
	\label{lem: Fraenkel_1}	
	Let $f\in C^{0,\alpha}(B_R)$. For any $x\in \mathrm{Int}(B_R)$,
	$$\partial_{ij}\mathcal{N}(f)(x)=\int_{B_R} \partial_{ij}\Gamma(x-y) \big ( f(y)-f(x) \big )dy-\frac{\delta_{ij}}{n}f(x),$$ where $\mathrm{Int}(B_R)$ is the interior of $B_R$.
\end{lemma}

\vspace{.1in}

Now we discuss a technical result in order to study the H\" older estimate for functions under the operator $\mathcal{N}$. Let $x$ be any interior point of $B_R$, and let $B_{\rho}(x)$ be the open ball centered at $x$ with radius $\rho$.  It was proved in  \cite{Fraenkel} that

\begin{lemma}  
	\label{lem: Fraenkel_2}
	There is a constant $C(n)$ depending only on $n$, such that for any $x \in \mathrm{Int}(B_R)$ and $1\leq i, j\leq n$,
	$$\bigg|\int_{B_R \setminus B_{\rho}(x)} \partial_{ij} \Gamma(x-y)dy\bigg|\leq C(n).$$
\end{lemma}

\vspace{.1in}

\noindent
Although in \cite{Fraenkel} Lemma \ref{lem: Fraenkel_2} was proved in the case $\rho \leq \frac{R}{4}$, the proof actually holds for all $\rho$ without restrictions.

\vspace{.2in}

It is well known (see \cite{GT}) that $\mathcal{N}$ maps $\displaystyle C^{0,\alpha}\left (B_{\frac{R}{2}} \right )$ to $C^{0,\alpha} (B_R)$ continuously. Next we will prove a stronger estimate which is essential to our construction of the contraction map.

\vspace{.1in}

\begin{thm} 
	\label{thm: N(f)_bounded_by_f}
	If $f\in C^{0,\alpha}(B_R)$, then $\mathcal{N}(f)\in C^{2,\alpha}(B_R)$ and there is a constant $C(n,\alpha)$, independent of $R$, such that
	$$\|\mathcal{N}(f)\|^{(2, \alpha)}\leq C(n,\alpha) \|f\|_{\alpha}.$$
\end{thm} 

\vspace{.1in}
\noindent
We would like to point out that although the constant $C(n, \alpha)$ is independent of $R$, there is an $R^{\alpha}$ term in the definition of the weighted norm $\|f\|_{\alpha}$, so the $\|\mathcal{N}(f)\|^{(2, \alpha)} $ norm actually does depend on $R$. This shows the advantage of choosing the weighted norm over unweighted norm.

\vspace{.1in}
\pf Recall that by (\ref{eqn:(k)-norm}),
$$\|\mathcal{N}(f)\|^{(2,\alpha)}= \max_{|\beta|=2}  \big \{\|\partial^\beta \mathcal{N}(f)\|_{\alpha} \big \} =  \max_{1 \leq i,j \leq n}  \big \{\|\partial _{ij} \mathcal{N}(f)\|_{\alpha} \big \} $$
Let
$$\zeta(x)=\int_{B_R} \partial_{ij}\Gamma(x-y)\big ( f(y)-f(x) \big )dy.$$
By Lemma \ref{lem: Fraenkel_1}, \,\,  $\displaystyle \partial_{ij}\mathcal{N}(f)=\zeta(x)-\frac{\delta_{ij}}{n}f(x)$, and therefore
\begin{equation}
	\label{eqn:partial_ij_N_bound}
	\|\partial_{ij}\mathcal{N}(f)\|_{\alpha} \leq \|\zeta\|_{\alpha}+\|f\|_{\alpha}.
\end{equation}

\vspace{.1in}
\noindent
So we only need to bound $\|\zeta\|_{\alpha}$ in terms of $\|f\|_{\alpha}$. 

\vspace{.1in}
\noindent
First, for $x\in \mathrm{Int} (B_R)$,

\allowdisplaybreaks
\begin{eqnarray}
	\label{eqn:phi_bound}
	|\zeta(x)|&\leq& \int_{B_R} |\partial_{ij}\Gamma(x-y)||f(y)-f(x)|dy\nonumber\\
		&\leq & CH_\alpha[f] \int_{B_R} \frac{|y-x|^\alpha}{|x-y|^n}dy\nonumber\\
	&\leq & CH_\alpha[f] \int_{B_{2R}(x)} \frac{|y-x|^\alpha}{|x-y|^n}dy\nonumber\\
	&= & CH_\alpha[f]\int_{S^{n-1}}\int_0^{2R}\frac{r^\alpha r^{n-1}}{r^n}drd\sigma \nonumber \\
	& = & C(n,\alpha)H_\alpha[f] R^\alpha \nonumber \\
	& \leq & C(n,\alpha) \|f\|_{\alpha},
\end{eqnarray}
where the polar coordinates are centered at $x$.

\vspace{.1in}
\noindent
To compute the H\"older constant of $\zeta$, let $x, x'$ be two (distinct) points in $B_R$. Let $B_{\rho}(x)$ be the open ball of radius $\rho=2|x-x'|$ and centered at $x$. 

\allowdisplaybreaks
\begin{eqnarray}
	\label{eqn:zeta_Holder}
	\zeta(x)-\zeta(x')&=&\int_{B_R} \partial_{ij}\Gamma(x-y) \Big ( f(y)-f(x) \Big )dy-\int_{B_R} \partial_{ij}\Gamma(x'-y) \Big ( f(y)-f(x') \Big )dy \nonumber\\
	&=&\int_{B_R\setminus B_{\rho}(x)} \partial_{ij}\Gamma(x-y)\Big ( f(y)-f(x) \Big )dy \nonumber \\
	& & +\int_{B_R \cap B_{\rho}(x)} \partial_{ij}\Gamma(x-y)\Big ( f(y)-f(x) \Big )dy \nonumber\\
	&-&\int_{B_R \setminus B_{\rho}(x)} \partial_{ij}\Gamma(x'-y) \Big ( f(y)-f(x') \Big )dy \nonumber \\
	& - & \int_{B_R \cap B_{\rho}(x)} \partial_{ij}\Gamma(x'-y) \Big (f(y)-f(x') \Big )dy\nonumber\\
	&=&\underbrace{ \int_{B_R\setminus B_{\rho}(x)} \Big(\partial_{ij}\Gamma(x-y)-\partial_{ij}\Gamma(x'-y) \Big )\Big ( f(y)-f(x) \Big )dy }_{I_1}  \\
	&+&\underbrace{ \Big (f(x')-f(x) \Big )\int_{B_R\setminus B_{\rho}(x)} \partial_{ij}\Gamma(x'-y)dy }_{I_2} \nonumber\\
	&+&\underbrace{ \int_{B_R \cap B_{\rho}(x)} \partial_{ij}\Gamma(x-y) \Big (f(y)-f(x) \Big )dy }_{I_3}  \nonumber \\
	& - &  \underbrace{ \int_{B_R \cap B_{\rho}(x)} \partial_{ij}\Gamma(x'-y)\Big (f(y)-f(x') \Big )dy }_{I_4}.  \nonumber
\end{eqnarray}

\vspace{.1in}

\noindent
Next we will estimate each of $I_1$, $I_2$, $I_3$, and $I_4$.

\allowdisplaybreaks
\begin{eqnarray*}
	|I_1|&=&\bigg|\int_{B_R\setminus B_{\rho}(x)} \Big (\partial_{ij}\Gamma(x-y)-\partial_{ij}\Gamma(x'-y) \Big )\Big (f(y)-f(x) \Big )dy\bigg|\nonumber\\
	&\leq & H_\alpha[f]|x-x'|\int_{B_R\setminus B_{\rho}(x)}\big |\nabla\partial_{ij}\Gamma(\hat{x}-y) \big ||y-x|^\alpha dy,
\end{eqnarray*}
where $\hat{x}$ is a point on the line segment between $x, x'$.

\vspace{.1in}
\noindent
Since $y \in B_R\setminus B_{\rho}(x)$, $$ |x-y| > \rho = 2 |x-x'| \geq 2|x-\hat{x}|,$$ then 
\begin{equation*}
	|\hat{x}-y| \geq |x-y| - |x-\hat{x} | \geq  \frac{1}{2} |x-y|$$ and $$ \big |\nabla\partial_{ij}\Gamma(\hat{x}-y) \big |  \leq  C(n) |\hat{x}-y|^{-n-1} \leq   C(n)  |x-y|^{-n-1}.
\end{equation*}

\noindent
Therefore,
 \allowdisplaybreaks
\begin{eqnarray}	
\label{eqn:I1}
|I_1|	&\leq& C(n) H_\alpha[f]|x-x'|\int_{B_R\setminus B_{\rho}(x)}|x-y|^{-n-1+\alpha}dy \nonumber \\
&\leq& C(n) H_\alpha[f]|x-x'|\int_{B_{2R}(x)\setminus B_{\rho}(x)}|x-y|^{-n-1+\alpha}dy \nonumber \\
	&=& C(n) H_\alpha[f]|x-x'|\int_{S^{n-1}}\int_\rho^{2R}  r^{-n-1+\alpha} r^{n-1}drd\sigma \nonumber \\
	& = & C(n, \alpha) H_\alpha[f]|x-x'| \big ( \rho ^{\alpha -1} - (2R)^{\alpha -1} \big ) \nonumber \\
	& < & C(n, \alpha) H_\alpha[f] |x-x'|^\alpha,
\end{eqnarray}

\noindent
where the polar coordinates are centered at $x$ and we have used $\rho=2|x-x'|$.

\vspace{.1in}
\noindent
By Lemma \ref{lem: Fraenkel_2}, we have
\begin{eqnarray}
	\label{eqn:I2}
	|I_2|&=&\bigg|\Big ( f(x')-f(x) \Big )\int_{B_R\setminus B_{\rho}(x)} \partial_{ij}\Gamma(x'-y)dy\bigg|\nonumber\\
	&\leq & C(n, \alpha) H_\alpha[f] |x-x'|^\alpha.
\end{eqnarray}

\vspace{.1in}
\noindent
The next term
\begin{eqnarray}
	\label{eqn:I3}
	|I_3|&=&\bigg|\int_{B_R \cap B_{\rho}(x)} \partial_{ij}\Gamma(x-y)\Big ( f(y)-f(x) \Big )dy\bigg| \nonumber\\
	&\leq & H_\alpha[f] \int_{B_R \cap B_{\rho}(x)}|\partial_{ij}\Gamma(x-y)||x-y|^\alpha dy\nonumber\nonumber\\
	&\leq & C(n)H_\alpha[f] \int_{B_R \cap B_{\rho}(x)}|x-y|^{-n}||x-y|^\alpha dy\nonumber\\
	&\leq& C(n) H_\alpha[f]\int_{S^{n-1}} \int_0^{\rho}  r^{-n+\alpha} r^{n-1}dr d\sigma \nonumber \\
	& \leq &C(n, \alpha)H_\alpha[f] |x-x'|^\alpha,
\end{eqnarray}
where the polar coordinates are also centered at $x$.

\vspace{.1in}
\noindent
The last term
\begin{eqnarray*}
	\label{eqn:I4}
	|I_4|&=&\bigg|\int_{B_R \cap B_{\rho}(x)} \partial_{ij}\Gamma(x'-y)\Big ( f(y)-f(x') \Big )dy\bigg| \nonumber \\
	&\leq & H_\alpha[f] \int_{B_R \cap B_{\rho}(x)}|\partial_{ij}\Gamma(x'-y)||x'-y|^\alpha dy \nonumber\\
	&\leq & C(n)H_\alpha[f] \int_{B_R \cap B_{\rho}(x)}|x'-y|^{-n}||x'-y|^\alpha dy.
\end{eqnarray*}

\noindent
For any $y \in B_R \cap B_{\rho}(x)$, $$ |y-x'| \leq |y-x|+|x-x'| \leq \rho + \frac{1}{2} \rho = \frac{3}{2}\rho.$$  

\noindent
Then
\begin{eqnarray*}
	\int_{B_R \cap B_{\rho}(x)}|x'-y|^{-n+\alpha} dy & \leq & \int_{ B_{\frac{3}{2}\rho}(x')}|x'-y|^{-n+\alpha} dy  \\
	& = & C(n)\int_0^{\frac{3}{2}\rho}  r^{-n+\alpha} r^{n-1}dr \\
	& = & C(n, \alpha)|x-x'|^\alpha, 
\end{eqnarray*}	
where the polar coordinates are centered at $x'$.

\vspace{.1in}
\noindent
Therefore, 
\begin{equation}
	\label{eqn:I4}
	|I_4| \leq C(n, \alpha) H_\alpha[f] |x-x'|^\alpha.
\end{equation}
Combining (\ref{eqn:zeta_Holder}) to (\ref{eqn:I4}), we have $$ H_{\alpha} [\zeta] \leq C(n, \alpha) H_\alpha[f]  .$$  This and (\ref{eqn:phi_bound}) now imply $$\|\zeta \|_{\alpha} = \sup_{x\in B_R}|\zeta(x)|+(2R)^\alpha H_\alpha[\zeta] \leq C(n, \alpha) \|f\|_{\alpha}, $$ 

\noindent
and thus by (\ref{eqn:partial_ij_N_bound}) the proof is completed.

\stop


\section{The Integral System and the Map}
\label{sec:Poisson}

\vspace{.2in}
A crucial observation in the proof of Theorem \ref{thm: main} is that we only need to prove it in the case $$\bm{c_0}=\bm{0} \hspace{.2in} \text{and} \hspace{.2in} \bm{c_1}=\bm{0}. $$
For arbitrary $\bm{c_0} \in B'_{R'}$ and $\bm{c_1} \in \mathbb{R}^{mn}$, we may choose $R$ small enough so that $\bm{c_0}+\bm{c_1}x \in B'_{R'}$ for all $x \in B_R$. We first solve the elliptic quasilinear system

\begin{equation*}
	\renewcommand{\arraystretch}{1.5}
	\left \{
	\begin{array}{r@{}l}
			\displaystyle \sum _{i,j=1}^na^{ij}\big ( x, \bm{v}(x)+\bm{c_0}+\bm{c_1}x, D\bm{v}(x)+\bm{c_1}\big ) D_{ij} v^k &{} = \phi^k\big ( x, \bm{v}(x)+\bm{c_0}+\bm{c_1}x, D\bm{v}(x)+\bm{c_1}\big ) \\
		\bm{v}(0) &{} = \bm{0} \\
		D\bm{v}(0) &{} =  \bm{0}
	\end{array}
	\right.
\end{equation*}
for $k=1,...,m$.  Then, $$ \bm{u}(x)=\bm{v}(x)+\bm{c_0}+\bm{c_1}x$$ will be a solution to (\ref{eqn:main}).

\vspace{.1in}
\noindent
Thus in the rest of the proof we will assume $$\bm{c_0}=\bm{0} \,\, \,\, \text{and} \,\, \,\, \bm{c_1}=\bm{0}.$$

\vspace{.2in}

Consider the elliptic system (\ref{eqn:main}), for any $1 \leq k \leq m$, we can write
$$  Lu^k= \sum _{i,j=1}^n a^{ij}\left ( 0, \bm{0}, \bm{0} \right ) D_{ij} u^k - \sum _{i,j=1}^n  \Big [ a^{ij}\left ( 0, \bm{0}, \bm{0} \right )  - a^{ij} \left ( x, \bm{u}(x), D\bm{u}(x)\right ) \Big ] D_{ij} u^k .  $$

\vspace{.1in}
\noindent
Therefore, (\ref{eqn:main}) can be written in vector form as 

\begin{eqnarray}
	\label{eqn:vector_constant}
	\sum_{i,j=1}^n a^{ij} \left ( 0, \bm{0}, \bm{0}  \right ) D_{ij} \bm{u} & = & \bm{\phi} \big (x, \bm{u}(x), D\bm{u}(x) \big ) \nonumber \\
	& + &   \displaystyle \sum _{i,j=1}^n  \Big [ a^{ij}\left ( 0, \bm{0}, \bm{0} \right )  - a^{ij} \left ( x, \bm{u}(x), D\bm{u}(x)\right ) \Big ] D_{ij} \bm{u}.
\end{eqnarray}

\noindent
Since the constant matrix $\displaystyle \Big ( a_{ij} \left ( 0, \bm{0}, \bm{0} \right ) \Big )_{n \times n}$ is positive definite, after a change of the coordinates  $\tilde{x}=Px$, where $P=(p_{ij})_{n \times n}$ is an invertible $n \times n$ matrix, we can write 
$$\sum_{i,j=1}^n a^{ij} \left ( 0, \bm{0}, \bm{0}  \right ) D_{ij} \bm{u} (x) = \displaystyle \sum _{i=1}^n \frac{\partial ^2 \bm{u}}{\partial \tilde{x}_i ^2}(\tilde{x}). $$  

\vspace{.1in}
\noindent
The right-hand side of (\ref{eqn:vector_constant}) becomes
\begin{eqnarray*}
& &   \bm{\phi} \big (P^{-1}\tilde{x}, \bm{u}(P^{-1}\tilde{x}), \tilde{D}\bm{u}(\tilde{x})P \big ) \\
& + &   \displaystyle \sum _{i,j=1}^n  \Big [ a^{ij}\left ( 0, \bm{0}, \bm{0} \right )  - a^{ij} \left ( P^{-1}\tilde{x}, \bm{u}(P^{-1}\tilde{x}), \tilde{D}\bm{u}(\tilde{x})P\right ) \Big ]    \displaystyle \sum _{k,l=1}^n p_{li}\frac{\partial ^2 \bm{u}}{\partial \tilde{x}_k \partial \tilde{x}_l}(\tilde{x}) p_{kj}, 
\end{eqnarray*}
where $\tilde{D}$ denotes differentiation with respect to the $\tilde{x}$ variable.

\vspace{.1in}
\noindent
After swapping the indices $(i,j)$ with $(k,l)$, it can be expressed as

\begin{eqnarray*}
&  &  \bm{\phi} \left (P^{-1}\tilde{x}, \bm{u}(P^{-1}\tilde{x}), \tilde{D}\bm{u}(\tilde{x})P \right )  \\
& + &   \displaystyle \sum _{i,j=1}^n \sum _{k,l=1}^n \Big [ a^{kl}\left ( 0, \bm{0}, \bm{0} \right )  - a^{kl} \left ( P^{-1}\tilde{x}, \bm{u}(P^{-1}\tilde{x}), \tilde{D}\bm{u}(\tilde{x})P\right ) \Big ]    p_{jk}p_{il} \frac{\partial ^2 \bm{u}}{\partial \tilde{x}_i \partial \tilde{x}_j}(\tilde{x}). 
 \end{eqnarray*}

\vspace{.1in}
\noindent
Let $$ \displaystyle \bm{\psi} \left ( \tilde{x}, \bm{u}(\tilde{x}), \tilde{D}\bm{u}(\tilde{x}) \right ) = \bm{\phi} \left ( P^{-1}\tilde{x}, \bm{u}(P^{-1}\tilde{x}), \tilde{D}\bm{u}(\tilde{x})P \right ) $$ and let 
$$\displaystyle b^{ij } \left (\tilde{x}, \bm{u}(\tilde{x}), \tilde{D}\bm{u}(\tilde{x}) \right ) = \displaystyle  \sum _{k,l=1}^n \Big [ a^{kl}\left ( 0, \bm{0}, \bm{0} \right )  - a^{kl} \left ( P^{-1}\tilde{x}, \bm{u}(P^{-1}\tilde{x}), \tilde{D}\bm{u}(\tilde{x})P\right ) \Big ]    p_{jk}p_{il}.  $$

\vspace{.1in}

\noindent
Then we can denote the right-hand side of (\ref{eqn:vector_constant}) as $$ \displaystyle \bm{\psi} \left ( \tilde{x}, \bm{u}(\tilde{x}), \tilde{D}\bm{u}(\tilde{x}) \right ) + \displaystyle \sum _{i,j=1}^n b^{ij } \left ( \tilde{x}, \bm{u}(\tilde{x}), \tilde{D}\bm{u}(\tilde{x}) \right )  \tilde{D}_{ij} (\bm{u}),  $$ where $\displaystyle \tilde{D}_{ij} (\bm{u}) = \frac{\partial ^2 \bm{u}}{\partial \tilde{x}_i \partial \tilde{x}_j}(\tilde{x}). $

\vspace{.1in}

\noindent
Thus in the $\tilde{x}$ coordinates (\ref{eqn:vector_constant}) becomes 
$$\displaystyle \sum _{i=1}^n \frac{\partial ^2 \bm{u}}{\partial \tilde{x}_i ^2}(\tilde{x}) = \displaystyle \bm{\psi} \left ( \tilde{x}, \bm{u}(\tilde{x}), \tilde{D}\bm{u}(\tilde{x}) \right ) + \displaystyle \sum _{i,j=1}^n b^{ij } \left ( \tilde{x}, \bm{u}(\tilde{x}), \tilde{D}\bm{u}(\tilde{x}) \right )  \tilde{D}_{ij} (\bm{u}).  $$

\noindent
Since the two coordinate systems have the same origin, in the $\tilde{x}$ coordinates we still have $$\bm{u}(0)=\bm{0} \hspace{.2in} \text{  and } \hspace{.2in} \tilde{D}\bm{u}(0)=\bm{0}\cdot P^{-1} = \bm{0}, $$ and thus
\begin{equation*}
	\label{eqn:b_ij_0}
	b^{ij}(0, \bm{0}, \bm{0})= 0.
\end{equation*} 

\vspace{.1in}

\noindent
Therefore, we only need to prove Theorem \ref{thm: main} under the coordinate system $\tilde{x}$, but to simplify notations we will drop the $\sim$ in all the subsequent notations.

\vspace{.1in}

\noindent
Our goal now is to prove the following result about a Poisson type system.
\begin{lemma}
	\label{thm: Poisson}
	 Assume the functions $\bm{\psi} (x, p, q)= \big (\psi ^1(x, p, q), ..., \psi ^m(x, p, q)\big )$  and $b^{ij}(x,p,q) $  are in $C_{loc}^{1,\alpha} \left ( B_R \times B'_{R'} \times  \mathbb{R}^{mn} \right ) $.  If $ b^{ij}(0, \bm{0}, \bm{0})= 0,  $ then the quasi-linear system
	\begin{equation}
		\label{eqn:main_poisson}
		\left \{
		\renewcommand{\arraystretch}{1.5}
		\begin{array}{r@{}l}
		\Delta \bm{u} &{} = \bm{\psi}\big (x, \bm{u}(x), D\bm{u}(x) \big ) +  \displaystyle \sum _{i,j=1}^n  b^{ij} \big ( x, \bm{u}(x), D\bm{u}(x) \big ) D_{ij}\bm{u}\\
			\bm{u}(0) &{} =  \bm{0} \\
			D\bm{u}(0) &{} =  \bm{0}
		\end{array}
	\right.
	\end{equation}
	
	\noindent
	has $C^{2,\alpha}$ solutions from $B_R$ to $B'_{R'} $ when $R$ is sufficiently small.
	
\end{lemma}

\vspace{.1in}
Denote 
\begin{equation}
	\label{eqn:Big_psi}
	\bm{\Psi} \big ( x, \bm{u}(x), D\bm{u}(x), D^2\bm{u}(x) \big ) =  -\bm{\psi}\big (x, \bm{u}(x), D\bm{u}(x) \big ) -  \displaystyle \sum _{i,j=1}^n  b^{ij} \big ( x, \bm{u}(x), D\bm{u}(x) \big ) D_{ij}\bm{u},
\end{equation}
so the equation in (\ref{eqn:main_poisson}) can be written as 
\begin{equation}
	\label{eqn: vector_psi}
	\Delta \bm{u} = - \bm{\Psi} \big ( x, \bm{u}(x), D\bm{u}(x), D^2\bm{u}(x) \big ).
\end{equation}

\vspace{.1in}
\noindent
A key observation is that if we can find a function $\bm{u}=(u^1,..., u^m)$ satisfying 
\begin{equation}
	\label{eqn:integral_Psi}
	\bm{u} = \bm{h}+ \int_{B_R} \Gamma(x-y) \bm{\Psi} \big ( y, \bm{u}(y), D\bm{u}(y), D^2\bm{u}(y) \big )dy,
\end{equation}

\noindent
where $\bm{h}=(h^1,...,h^m)$  and each $h^i(x)$ is a harmonic function, then $\bm{u}$ will be a solution of (\ref{eqn: vector_psi}). Furthermore, because $\bm{h}$ is arbitrary, it allows us to construct infinitely many such solutions.

\vspace{.1in}
\noindent
 Written as a system of equations, (\ref{eqn:integral_Psi}) is
\begin{equation*}
	\label{eqn:integral_system}
	\left \{
		\renewcommand{\arraystretch}{2}
	\begin{array}{r@{}l}
		u^1 &{} =   h^1 + \displaystyle \int_{B_R} \Gamma(x-y) \Psi ^1 \big ( y, \bm{u}(y), D\bm{u}(y), D^2\bm{u}(y)  \big )dy,\\
		&{} \vdots  \\
			u^m &{}=   h^m + \displaystyle \int_{B_R} \Gamma(x-y) \Psi ^m \big ( y, \bm{u}(y), D\bm{u}(y), D^2\bm{u}(y)  \big )dy.
		\end{array}
			\right.			
	\end{equation*} 

\vspace{.1in}
Since our strategy is to construct a contraction map on functions in $C_0^{2, \alpha}$ and apply the Banach Fixed Point Theorem, it seems natural to define the map as the integral on the right hand side of (\ref{eqn:integral_Psi}).  However, the function may not remain in $C_0^{2, \alpha}$ under such a map,  so we need to modify it slightly.

\vspace{.1in}
\noindent
For any function $\bm{f}=(f^1,...,f^m) \in  C_0^{2, \alpha} (B_R) \times \cdots \times C_0^{2, \alpha} (B_R)  $, define 
\begin{eqnarray}
	\label{defn:omega}
	\omega^i(\bm{f})(x) & = & \mathcal{N}\Big (\Psi^i \big (y, \bm{f}(y), D\bm{f}(y), D^2\bm{f}(y)\big ) \Big )(x)   \\
	& = &  \int_{B_R} \Gamma(x-y) \Psi ^i \big ( y, \bm{f}(y), D\bm{f}(y), D^2\bm{f}(y) \big )dy \nonumber
\end{eqnarray}
and 
\begin{equation}
	\label{eqn:theta}
	\Theta^i(\bm{f})(x)=\omega^i(\bm{f})(x)-\omega^i(\bm{f})(0)-\sum_{j=1}^n\partial_j\big (\omega^i(\bm{f})\big )(0)x_j-\frac{1}{2}\sum_{\substack{k,l=1\\k\not =l }}^n\partial_k\partial_l\big (\omega^i(\bm{f})\big )(0)x_kx_l.
\end{equation}

\vspace{.1in}
\noindent
Since $\Psi^i $ is $  C^{0, \alpha} $, by Theorem \ref{thm: N(f)_bounded_by_f} we know $\omega^i(\bm{f}) $ is $  C^{2, \alpha},$ and hence $\Theta^i(\bm{f}) $ is $  C^{2, \alpha}.$

\vspace{.1in}
\noindent
Direct calculations of the first and second derivatives show that  
  $$  \displaystyle  \Theta^i(\bm{f}) \in C_0^{2, \alpha} (B_R) $$ and $$ \Delta \Theta^i(\bm{f}) = \Delta \omega^i(\bm{f})  = - \Psi ^i \big ( x, \bm{f}(x), D\bm{f}(x), D^2\bm{f}(x) \big )  .$$

\vspace{.1in}
\noindent
Thus, if a function $\bm{u}$ satisfies $$ \Theta^i(\bm{u})=\bm{u}$$ for all $i=1,...,m,$ then it will be a solution to (\ref{eqn: vector_psi}).  Therefore, we define the main operator $$\bm{\Theta}: C_0^{2, \alpha} (B_R) \times \cdots \times C_0^{2, \alpha} (B_R) \to  C_0^{2, \alpha} (B_R) \times \cdots \times C_0^{2, \alpha} (B_R) $$ by
$$ \bm{\Theta}(\bm{f}) = \big ( \Theta ^1(\bm{f}), ..., \Theta ^m(\bm{f}) \big ). $$ 

\vspace{.1in}
\noindent
We will find a solution to the system (\ref{eqn: vector_psi}) by applying the Banach Fixed Point Theorem to the map $\bm{\Theta}$ on a closed subset of $C_0^{2, \alpha} (B_R) \times \cdots \times C_0^{2, \alpha} (B_R) $, defined by
$$ \bm{\mathcal{E}}(R,\gamma)= \left \{ \bm{f} \in C_0^{2, \alpha} (B_R) \times \cdots \times C_0^{2, \alpha} (B_R) \Big |  \| \bm{f} \|^{(2, \alpha) }   \leq \gamma \right \},    $$ where $\gamma >0$ is a parameter to be determined.

\vspace{.1in}
\noindent
In order to show that $\bm{\Theta}$ maps $ \bm{\mathcal{E}}(R,\gamma)$ into itself, we need to estimate $\|\bm{\Theta}(\bm{f} )\|^{(2, \alpha)}$.  To show $\bm{\Theta}$ is a contraction, we need to estimate  $  \|\bm{\Theta}(\bm{f}-\bm{g})\|^{(2, \alpha)} $.  The essential step in the estimates is to choose an appropriate value for $\gamma$ and then adjust $R$ accordingly.
 \vspace{.2in}
 
 \noindent
 In the rest of the paper, we will use $C$ to denote all constants that depend only on $m,n,$ and $\alpha$.

  \vspace{.2in}
  
\section{The Estimate for $\|\bm{\Theta}(\bm{f} )\|^{(2, \alpha)}$}
\label{sec:estimate_f}

\vspace{.2in}

First, we show that $\bm{\Theta}$ maps $ \bm{\mathcal{E}}(R,\gamma)$ into itself when $\gamma$ is appropriately chosen.  Specifically,  we will show there is a $\gamma >0$ such that for all $\bm{f}$ in $ C_0^{2, \alpha} (B_R) \times \cdots \times C_0^{2, \alpha} (B_R)$ with $\| \bm{f} \|^{(2, \alpha) } \leq \gamma$, we have $\displaystyle \| \Theta ^i (\bm{f}) \|^{(2, \alpha) } \leq  \frac{\gamma}{2}$.  Then, by the definition of the norm of vector functions in (\ref{defn:vector_norm}) we will have $\displaystyle \|\bm{\Theta}(\bm{f} )\|^{(2, \alpha)} \leq \frac{\gamma}{2}. $

\vspace{.1in}
\noindent
From (\ref{eqn:theta}) we have 
\begin{equation}
	\label{eqn:theta_norm}
	\|\Theta^i(\bm{f})\|^{(2, \alpha)}\leq \|\omega^i(\bm{f})\|^{(2,\alpha)}+\frac{1}{2}\Bigg | \sum_{\substack{k,l=1\\k\not =l }}^n\partial_k\partial_l\big (\omega^i(\bm{f})\big )(0) \Bigg |.
\end{equation}

\noindent
It follows from (\ref{defn:omega}) and Theorem \ref{thm: N(f)_bounded_by_f} that

\begin{equation}
	\label{eqn: omega_estimate}
	\|\omega^i(\bm{f})\|^{(2,\alpha)} \leq C \|\Psi^i \big (x, \bm{f}(x), D\bm{f}(x), D^2\bm{f}(x)\big ) \|_{\alpha}.
\end{equation}

\vspace{.1in}
\noindent
Also by Theorem \ref{thm: N(f)_bounded_by_f}, we know
\begin{eqnarray}
	\label{eqn:omega_0_estimate}
\Big | \partial_k\partial_l\big (\omega^i(\bm{f})\big )(0) \Big | & \leq &   \sup_{B_R}\Big | \partial_k\partial_l\big (\omega^i(\bm{f})\big ) \Big | \nonumber \\
& \leq & \|\omega^i(\bm{f})\|^{(2,\alpha)} \nonumber \\
& \leq &  C \|\Psi^i \big (x,\bm{f}(x), D\bm{f}(x), D^2\bm{f}(x)\big ) \|_{\alpha}. 
\end{eqnarray}

\vspace{.1in}
\noindent
Therefore we only need to estimate $\displaystyle  \|\Psi^i \big (x, \bm{f}(x), D\bm{f}(x), D^2\bm{f}(x)\big ) \|_{\alpha}.  $

\vspace{.1in}
\noindent
We use the coordinates $x=(x_1,...,x_n),$ $p=(p_1,...,p_m),$ $\displaystyle q=\left (q_k^j \right)_{\substack{j=1,...,m\\k=1,...,n}},$ and $\displaystyle  r=\left (r_{kl}^j \right )_{\substack{j=1,...,m\\k,l=1,...,n }}$ to denote the components in $\Psi ^i (x,p,q,r) $.
\vspace{.1in}
\noindent
Then

\allowdisplaybreaks
\begin{eqnarray}
	\label{eqn:psi_estimate_1}
	&&\Psi^i \Big (x,\bm{f}(x), D\bm{f}(x), D^2\bm{f}(x) \Big ) - \Psi^i(0,\bm{0},...,\bm{0})\nonumber\\
	&=&\int_0^1 \frac{d}{dt}\Psi^i \Big ( tx,t\bm{f}(x),tD\bm{f}(x), tD^2\bm{f}(x) \Big )dt\nonumber\\
	&=& \sum_{j=1}^n E_j x_j+\sum_{j=1}^m A_jf^j+\sum_{k=1}^n\sum_{j=1}^m B_k^j\partial_k f^j + \sum_{j=1}^m \sum_{k,l=1}^nC_{kl}^j\partial_k\partial_l f^j,
\end{eqnarray}
where
\allowdisplaybreaks
\begin{equation}
	\label{eqn:const_deriv_psi_A}
	\begin{aligned}
			E_j &  = & \displaystyle \int_0^1 \frac{\partial}{\partial x_j}\Psi^i(tx, t \bm{f}, tD\bm{f}, tD^2\bm{f})dt,\\
	A_j & = & \displaystyle \int_0^1 \frac{\partial}{\partial p_j}\Psi^i(tx, t \bm{f}, tD\bm{f} , tD^2\bm{f})dt,\\
	\smallskip
	B_k^j&  = & \displaystyle \int_0^1 \frac{\partial}{\partial q_k^j}\Psi^i(tx, t \bm{f}, tD\bm{f}, tD^2\bm{f})dt,\\
	C_{kl}^j& = & \displaystyle \int_0^1 \frac{\partial}{\partial r_{kl}^j}\Psi^i(tx, t \bm{f}, tD\bm{f}, tD^2\bm{f})dt.
\end{aligned}
\end{equation}

\vspace{.1in}
\noindent
By Lemma \ref{lem:f-bounded-by_f(k)}, Lemma \ref{lem: compare_norms}, and the fact that $\| \bm{f} \| ^{2, \alpha} 
\leq \gamma$, 

\allowdisplaybreaks
\begin{equation}
	\label{eqn:consant_domain}
	\addtolength{\jot}{1em}
	\begin{aligned}
		\|f^j\|_{\alpha} \leq {}& CR^2 \|f^j\|^{(2, \alpha)} \leq CR^2\gamma,  \\
		\|\partial _{k} f^j\|_{\alpha} \leq   {}& \|f^j\|^{(1, \alpha)} \leq CR \|f^j \|^{(2, \alpha)} \leq CR\gamma,  \\
		\| \partial_k\partial_l f^j \|_{\alpha} \leq  {}& \|f^j\|^{(2, \alpha)}  \leq \gamma.
	\end{aligned}
\end{equation}


\vspace{.2in}
\noindent
Therefore, by (\ref{eqn:psi_estimate_1}), (\ref{eqn:Banach_algebra}), and the fact $\displaystyle \|  x_j \|_{\alpha} =3R $, we have 

\allowdisplaybreaks
\begin{eqnarray}
	\label{eqn:psi_estimate_2}
	  \Big \|\Psi^i \Big (x,\bm{f}(x), D\bm{f}(x), D^2\bm{f}(x) \Big ) \Big \|_{\alpha}  
	& \leq & \big | \Psi^i  (0,\bm{0}, \bm{0}, \bm{0}  ) \big | + 3 R\sum_{j=1}^n\|E_j\|_{\alpha}  + CR^2  \gamma\sum_{j=1}^m\|A_j\|_{\alpha} \nonumber \\
	& + & CR \gamma \sum_{k=1}^n\sum_{j=1}^m \|B_k^j\|_{\alpha}  + \gamma\sum_{j=1}^m \sum_{k,l=1}^n\|C^j_{kl}\|_{\alpha}.  
\end{eqnarray}

\vspace{.2in}
\noindent
Next, we will estimate $ \|A_j\|_{\alpha} $, $\|B_k^j\|_{\alpha} $,  $\|C^j_{kl}\|_{\alpha} $, and $ \|E_j\|_{\alpha} $.  
\
By (\ref{eqn:consant_domain}) we only need to estimate on the domain
\begin{equation}
	\label{eqn:constant_domain}
	\mathcal{B}(R, \gamma)  =  \Big \{ (x,p,q,r) \big | x \in B_R, |p|\leq CR^2\gamma,|q|\leq CR\gamma, |r|\leq \gamma \Big  \}.
\end{equation}

\vspace{.1in}

\noindent
Denote
\allowdisplaybreaks
\begin{equation}
	\label{eqn:constants_bound_A}
	\begin{aligned} 
	A[R,\gamma]= {}&\max \left \{\bigg|\frac{\partial \Psi^i}{\partial p_j}\bigg|_{\mathcal{B}(R,\gamma)}, \,\,i=1,..., m; \,\, j=1,...,m \right \}, \\
	H_\alpha^A[R,\gamma] ={} &\max \left \{H_\alpha\bigg[\frac{\partial \Psi^i}{\partial p_j}\bigg]_{\mathcal{B}(R,\gamma)}, \,\, i=1,..., m; \,\, j=1,...,m \right \},
		\end{aligned}
\end{equation}

\allowdisplaybreaks
\begin{equation}
	\label{eqn:constants_bound_B}
	\begin{aligned} 
	B[R,\gamma] = {} &\max \left  \{\bigg|\frac{\partial \Psi^i}{\partial q_k^j}\bigg|_{\mathcal{B}(R,\gamma)}, \,\,i,j=1,..., m; \,\, k=1,...,n \right \}, \\
	H_\alpha^B[R,\gamma] ={} &\max \left \{H_\alpha\bigg[\frac{\partial \Psi^i}{\partial q_k^j}\bigg]_{\mathcal{B}(R,\gamma)}, \,\,i,j=1,..., m; \,\, k=1,...,n \right \},
		\end{aligned}
\end{equation}

\allowdisplaybreaks
\begin{equation}
	\label{eqn:constants_bound_C}
	\begin{aligned} 
	C[R,\gamma] ={} &\max \left \{\bigg|\frac{\partial \Psi^i}{\partial r_{kl}^j}\bigg|_{\mathcal{B}(R,\gamma)}, \,\,i,j=1,..., m;\,\, k,l=1,...,n \right \}, \\
	H_\alpha^C[R,\gamma]= {} &\max \left \{ H_\alpha\bigg[\frac{\partial \Psi^i}{\partial r_{kl}^j}\bigg]_{\mathcal{B}(R,\gamma)}, \,\,i,j=1,..., m; \,\, k,l=1,...,n \right \},
		\end{aligned}
\end{equation}

\allowdisplaybreaks
\begin{equation}
	\label{eqn:constants_bound_E}
	\begin{aligned} 
	E[R,\gamma]= {}&\max \left \{\bigg|\frac{\partial \Psi^i}{\partial x_j}\bigg|_{\mathcal{B}(R,\gamma)}, \,\,i=1,..., m; \,\, j=1,...,n \right \}, \\
	H_\alpha^E[R,\gamma] ={} &\max \left \{H_\alpha\bigg[\frac{\partial \Psi^i}{\partial x_j}\bigg]_{\mathcal{B}(R,\gamma)}, \,\, i=1,..., m; \,\, j=1,...,n \right \}.
	\end{aligned}
\end{equation}

\vspace{.1in}
\noindent


\noindent
We use $H_1$ to denote the Lipschitz constant in the $r$ variable and define

\allowdisplaybreaks
\begin{equation}
\begin{aligned}
H_1^A[R,\gamma] = {}&\max \left \{ H_1\left [\frac{\partial \Psi ^i}{\partial p_j}\right ]_{\mathcal{B}(R,\gamma)}, \,\, i=1,..., m; \,, j=1,...,m \right \}, \\
H_1^B[R,\gamma] = {} &\max \left \{H_1\left [\frac{\partial \Psi^i}{\partial q_k^j}\right ]_{\mathcal{B}R,\gamma)}, \,\, i,j=1,..., m; \,\, k=1,...,n \right \}, \\
H_1^C[R,\gamma] = {}&\max \left \{H_1\left [\frac{\partial \Psi^i}{\partial r_{kl}^j}\right ]_{\mathcal{B}(R,\gamma)}, \,\, i,j=1,..., m; \,\, k,l=1,...,n \right \}, \\
H_1^E[R,\gamma] = {}&\max \left \{ H_1\left [\frac{\partial \Psi ^i}{\partial x_j}\right ]_{\mathcal{B}(R,\gamma)}, \,\, i=1,..., m; \,\, j=1,...,n \right \}.
\end{aligned}
\end{equation}

Note that by the definition of $ \bm{\Psi}$ in (\ref{eqn:Big_psi}), there is no $r$ variable in $\displaystyle \frac{\partial \Psi^i}{\partial r_{kl}^j}$, and therefore 

\begin{equation}
	\label{eqn:H_1_C=0}
	H_1^C[R,\gamma] =0.
\end{equation}

\vspace{.2in}
\subsection{Estimates for $\|A_j\|_{\alpha} $, $\|B_k^j\|_{\alpha} $,  $\|C^j_{kl}\|_{\alpha} $, and $\|E_j\|_{\alpha} $}

$ $

\vspace{.1in}
\noindent
By (\ref{eqn:const_deriv_psi_A}) and (\ref{eqn:constants_bound_A}), $$ \sup_{B_R} |A_j| \leq A[R, \gamma]. $$

Now we will estimate $H_\alpha[A_j]$. Let $x,x' \in B_R$,

\allowdisplaybreaks
\begin{eqnarray}
& & 	A_j(x)-A_j(x') \nonumber \\
& = &
	\displaystyle \int_0^1 \frac{\partial}{\partial p_j}\Psi^i(tx, t \bm{f}(x), tD\bm{f}(x) , tD^2\bm{f}(x))dt 
	- \displaystyle \int_0^1 \frac{\partial}{\partial p_j}\Psi^i(tx', t \bm{f}(x'), tD\bm{f}(x') , tD^2\bm{f}(x') )dt \nonumber\\
	& = &\displaystyle \int_0^1 \frac{\partial}{\partial p_j}\Psi^i(tx, t \bm{f}(x), tD\bm{f}(x) , tD^2\bm{f}(x))dt  - \displaystyle \int_0^1 \frac{\partial}{\partial p_j}\Psi^i(tx', t \bm{f}(x), tD\bm{f}(x) , tD^2\bm{f}(x))dt \nonumber\\
	& + &\displaystyle \int_0^1 \frac{\partial}{\partial p_j}\Psi^i(tx', t \bm{f}(x), tD\bm{f}(x) , tD^2\bm{f}(x))dt  -  \displaystyle \int_0^1 \frac{\partial}{\partial p_j}\Psi^i(tx', t \bm{f}(x'), tD\bm{f}(x) , tD^2\bm{f}(x))dt \nonumber\\
	& + &\displaystyle \int_0^1 \frac{\partial}{\partial p_j}\Psi^i(tx', t \bm{f}(x'), tD\bm{f}(x) , tD^2\bm{f}(x))dt  - \displaystyle \int_0^1 \frac{\partial}{\partial p_j}\Psi^i(tx', t \bm{f}(x'), tD\bm{f}(x') , tD^2\bm{f}(x))dt \nonumber \\
	& + & \displaystyle \int_0^1 \frac{\partial}{\partial p_j}\Psi^i(tx', t \bm{f}(x'), tD\bm{f}(x') , tD^2\bm{f}(x))dt  -  \displaystyle \int_0^1 \frac{\partial}{\partial p_j}\Psi^i(tx', t \bm{f}(x'), tD\bm{f}(x') , tD^2\bm{f}(x'))dt. \nonumber
\end{eqnarray}

\vspace{.1in}
\noindent
It then follows that
\allowdisplaybreaks
\begin{eqnarray*}
	&&|A_j(x)-A_j(x')| \\
	&\leq & H_\alpha^A \left [ R,\gamma \right ]|x-x'|^\alpha+H_\alpha^A \left [ R,\gamma \right ]\sum_{j=1}^m|f^j(x)-f^j(x')|^{\alpha}  \\
	&+&H_\alpha^A \left [R,\gamma \right ]\sum_{j=1}^m |D f^j(x)-Df^j(x')|^\alpha + H_1^A \left [R,\gamma \right ]\sum_{j=1}^m|D^2f^j(x)-D^2 f^j(x')| . 
\end{eqnarray*}

\noindent
Note that 
\begin{eqnarray*}
	|f^j(x)-f^j(x')| & \leq & \sup_{B_R}|Df^j| |x-x'| \\
	& \leq & \|f^j\|^{(1,\alpha)}|x-x'|\\
	& \leq & C R\|f^j\|^{(2,\alpha)}|x-x'|  \hspace{.2in} \text{ (by Lemma \ref{lem: compare_norms})}\\
	& \leq & C R \|\bm{f}\|^{(2,\alpha)}|x-x'| \\
	& \leq & C R\gamma |x-x'|.
\end{eqnarray*}
Similarly,
\begin{eqnarray*}
|D f^j(x)-Df^j(x')| & \leq & \sup_{B_R}|D^2f^j| |x-x'| \\
& \leq &  C \|\bm{f}\|^{(2,\alpha)}|x-x'|\\
 & \leq &  C \gamma |x-x'|.	
\end{eqnarray*}

\noindent 
In addition, 
\begin{eqnarray*}
	|D^2f^j(x)-D^2 f^j(x')| & \leq & H_{\alpha}\left [ D^2f^j \right ] |x-x'|^{\alpha} \\
	& \leq & (2R)^{-\alpha}\|f^j\|^{(2,\alpha)}   |x-x'|^{\alpha}  \hspace{.2in} \text{(by (\ref{eqn:(k)-norm}) )} \\
	& \leq & (2R)^{-\alpha}\|\bm{f}\|^{(2,\alpha)}   |x-x'|^{\alpha}\\
	& \leq & (2R)^{-\alpha}\gamma   |x-x'|^{\alpha}.
\end{eqnarray*}

\noindent
Therefore 
\begin{eqnarray*}
	&&|A_j(x)-A_j(x')| \\
	&\leq & H_\alpha^A \left [R,\gamma \right ]|x-x'|^\alpha + CH_\alpha^A \left [R,\gamma \right ]\cdot R^{\alpha}\gamma^{\alpha}  |x-x'|^{\alpha}   \\
	&+ &  CH_\alpha^A \left [R,\gamma \right ]\cdot \gamma^{\alpha}  |x-x'|^{\alpha} + CH_1^A \left [R,\gamma \right]\cdot  (2R)^{-\alpha}\gamma|x-x'|^{\alpha}, 
\end{eqnarray*}
which implies
\begin{eqnarray*}
	H_{\alpha}[A_j] & \leq & H_\alpha^A[R,\gamma] + CH_\alpha^A[R,\gamma]\cdot R^{\alpha}\gamma^{\alpha}  +  CH_\alpha^A[R,\gamma]\cdot\gamma^{\alpha}  +CH_1^A[R,\gamma]\cdot  (2R)^{-\alpha}\gamma.
\end{eqnarray*}

\noindent
Thus,
\begin{eqnarray*}
	\|A_j\|_{\alpha} & =& \sup_{B_R}|A_j|+(2R)^\alpha H_\alpha[A_j]   \\
	& \leq  & A[R,\gamma]+(2R)^\alpha \Big (H_\alpha^A[R,\gamma] + CH_\alpha^A[R,\gamma]\cdot R^{\alpha}\gamma^{\alpha}   
	 +  CH_\alpha^A[R,\gamma]\cdot\gamma^{\alpha}  \Big ) +  C\gamma H_1^A[R,\gamma ] .
\end{eqnarray*}

\noindent
Denote
\begin{eqnarray}
	\label{delta_B_R_gamma}
& & \delta _A(R,\gamma) \\
& = &  A[R,\gamma]+(2R)^\alpha \Big (H_\alpha^A[R,\gamma] + CH_\alpha^A[R,\gamma]\cdot R^{\alpha}\gamma^{\alpha}   
+  CH_\alpha^A[R,\gamma]\cdot\gamma^{\alpha}  \Big ) +  C\gamma H_1^A[R,\gamma ], \nonumber
\end{eqnarray}
then we can write
\begin{equation}
	\label{eqn:A_j_estimate}
		\|A_j\|_{\alpha} \leq \delta _A(R,\gamma). 
\end{equation}


\vspace{.2in}
\noindent
By similar calculations,

\begin{equation}
	\label{eqn:B_j_estimate}
	\|B^j_k\|_{\alpha} \leq \delta _B(R,\gamma),
\end{equation}
where
\begin{eqnarray}
	\label{eqn:delta_B_R_gamma}
	& & \delta _B(R,\gamma) \\
	 & = & B[R,\gamma]+(2R)^\alpha \Big (H_\alpha^B[R,\gamma] + CH_\alpha^B[R,\gamma]\cdot R^{\alpha}\gamma^{\alpha}   +  CH_\alpha^B[R,\gamma]\cdot\gamma^{\alpha}  \Big )  +  C\gamma H_1^B[R,\gamma] \nonumber;
\end{eqnarray}

\noindent
and
\begin{equation}
	\label{eqn:C_j_estimate}
	\|C^j_{kl}\|_{\alpha} \leq \delta _C(R,\gamma),
\end{equation}
where
\begin{eqnarray}
	\label{eqn:delta_C_R_gamma}
	& & \delta _C(R,\gamma) \\
	& = & C[R,\gamma]+(2R)^\alpha \Big (H_\alpha^C[R,\gamma] + CH_\alpha^C[R,\gamma]\cdot R^{\alpha}\gamma^{\alpha}  + CH_\alpha^C[R,\gamma]\cdot\gamma^{\alpha}  \Big )  +  C\gamma H_1^C[R,\gamma] \nonumber \\
	& = & C[R,\gamma]+(2R)^\alpha \Big (H_\alpha^C[R,\gamma] + CH_\alpha^C[R,\gamma]\cdot R^{\alpha}\gamma^{\alpha}  + CH_\alpha^C[R,\gamma]\cdot\gamma^{\alpha}  \Big ) \hspace{.3in} (\text{by \ref {eqn:H_1_C=0}}); \nonumber
\end{eqnarray}

\noindent
and we also have
\begin{equation}
	\label{eqn:E_j_estimate}
	\|E_j\|_{\alpha} \leq \delta _E(R,\gamma),
\end{equation}
where
\begin{eqnarray}
	\label{eqn:delta_E_R_gamma}
	& & \delta _E(R,\gamma) \\
	& = & E[R,\gamma]+(2R)^\alpha \Big (H_\alpha^E[R,\gamma] + CH_\alpha^E[R,\gamma]\cdot R^{\alpha}\gamma^{\alpha}   +  CH_\alpha^E[R,\gamma]\cdot\gamma^{\alpha}  \Big ) + C\gamma H_1^E[R,\gamma] \nonumber.
\end{eqnarray}

\vspace{.2in}

\subsection {The Estimates for $ \| \bm{\Theta} (\bm{f}) \|^{(2, \alpha)}$}

$ $

\vspace{.1in}
\noindent
Therefore, by (\ref{eqn:psi_estimate_2}),

\begin{eqnarray}
	\label{eqn:psi_estimate_3}
& &	\Big \|\Psi^i \Big (x,\bm{f}(x), \bm{Df}(x), \bm{D^2f}(x) \Big ) \Big \|_{\alpha} \\ 
	& \leq & \big | \Psi^i  (0,\bm{0}, \bm{0}, \bm{0}  ) \big | +  CR \delta _E(R,\gamma) + CR^2\gamma  \delta _A(R,\gamma)  +   C R\gamma \delta _B(R,\gamma) + C \gamma \delta _C(R,\gamma). \nonumber
\end{eqnarray}

\vspace{.1in} 
\noindent
Then by (\ref{eqn:theta_norm}), (\ref{eqn: omega_estimate}), (\ref{eqn:omega_0_estimate}), and (\ref{eqn:psi_estimate_3}), 
	
\begin{eqnarray}
		\|\Theta^i(\bm{f})\|^{(2, \alpha)} & \leq &  C \big | \Psi^i  (0,\bm{0}, \bm{0}, \bm{0}  ) \big | 
		+  C\Big (R \delta _E(R,\gamma) + R^2\gamma  \delta _A(R,\gamma) +   R\gamma \delta _B(R,\gamma) \Big ) + C\gamma \delta _C(R,\gamma). \nonumber 	
	\end{eqnarray}

\vspace{.1in} 
\noindent
By the expression of $\delta_C(R,\gamma)$ in (\ref{eqn:delta_C_R_gamma}), it can be written as

\allowdisplaybreaks
\begin{eqnarray*}
			\|\Theta^i(\bm{f})\|^{(2, \alpha)} & \leq &   C \big | \Psi^i  (0,\bm{0}, \bm{0}, \bm{0}  ) \big | +  C\Big (R \delta _E(R,\gamma) + R^2\gamma  \delta _A(R,\gamma) +   R\gamma \delta _B(R,\gamma) \Big )  + 	C\gamma C[R,\gamma]  \nonumber \\
& + & C\gamma (2R)^\alpha \Big (H_\alpha^C[R,\gamma] + CH_\alpha^C[R,\gamma]\cdot R^{\alpha}\gamma^{\alpha}  + CH_\alpha^C[R,\gamma]\cdot\gamma^{\alpha}  \Big ).		  
\end{eqnarray*}

\noindent
Denote 
\allowdisplaybreaks
\begin{eqnarray*}
	\epsilon(R,\gamma) & = & C\Big (R \delta _E(R,\gamma) + R^2\gamma  \delta _A(R,\gamma) +   R\gamma \delta _B(R,\gamma) \Big )  \\ 
	& + &C\gamma (2R)^\alpha \Big (H_\alpha^C[R,\gamma] + CH_\alpha^C[R,\gamma]\cdot R^{\alpha}\gamma^{\alpha}  + CH_\alpha^C[R,\gamma]\cdot\gamma^{\alpha}  \Big ).  
\end{eqnarray*}

\noindent
Then
\allowdisplaybreaks
\begin{equation}
	\label{eqn:theta_norm_bound_1}
	\|\Theta^i(\bm{f})\|^{(2, \alpha)}  \leq   C \big | \Psi^i  (0,\bm{0}, \bm{0}, \bm{0}  ) \big | + C\gamma C[R,\gamma] + 	\epsilon(R,\gamma),
\end{equation}
and
\allowdisplaybreaks
\begin{equation}
		\label{eqn:epsilon_R_gamma}
		\epsilon(R,\gamma) \to 0 \,\, \text{as} \,\, R\to 0.
\end{equation}

\vspace{.2in}
\noindent
Next, we will derive an upper bound for $C[R,\gamma]$.

\vspace{.2in}
\noindent
By (\ref{eqn:Big_psi}) 
$$\frac{\partial \Psi^i}{\partial r_{kl}^j}=
\begin{cases}
	b^{kl}\big (x, p, q \big ) & \text{if} \,\, i = j \\
	0 & \text{if} \,\, i \neq j.
\end{cases} 
 $$
So as defined in (\ref{eqn:constants_bound_C}), 
$$ C[R, \gamma])= \max \left \{\big|	b^{kl} (x, p, q)\big|_{\mathcal{B}(R,\gamma)} \right \}. $$

\vspace{.1in}
\noindent
Since $	b^{kl} (0, \bm{0}, \bm{0})=0,$ for any $(x,p,q,r)\in \mathcal{B}(R,\gamma)$ as defined in (\ref{eqn:constant_domain}),

\allowdisplaybreaks
\begin{eqnarray*}
\big|	b^{kl} (x, p, q)\big| & \leq & 	 \big|	b^{kl} (x, p, q) - b^{kl} (0, p, q) \big|  +   \big|	b^{kl} (0, p, q) - b^{kl} (0, \bm{0}, q) \big| + \big|	b^{kl} (0, \bm{0}, q) - b^{kl} (0, \bm{0}, \bm{0}) \big| \\
& \leq &  \sup_{\mathcal{B}(R,\gamma)} \bigg |  \frac{\partial b^{kl}}{\partial x}\bigg | \cdot R  + \sup_{\mathcal{B}(R,\gamma)} \bigg |  \frac{\partial b^{kl}}{\partial p}\bigg | \cdot C R^2\gamma  + \sup_{\mathcal{B}(R,\gamma)} \bigg |  \frac{\partial b^{kl}}{\partial q}\bigg | \cdot C R\gamma .
\end{eqnarray*}

\vspace{.1in}
\noindent
Therefore 
\begin{equation*}
	\label{eqn:C_R_gamma_bound}
	C[R, \gamma] \leq   \sup_{\mathcal{B}(R,\gamma)} \bigg |  \frac{\partial b^{kl}}{\partial x}\bigg | \cdot R   +  \sup_{\mathcal{B}(R,\gamma)} \bigg |  \frac{\partial b^{kl}}{\partial p}\bigg | \cdot C R^2\gamma  + \sup_{\mathcal{B}(R,\gamma)} \bigg |  \frac{\partial b^{kl}}{\partial q}\bigg | \cdot C R\gamma. 
\end{equation*}

\vspace{.1in}

\noindent
Thus we know 
\begin{equation}
	\label{eqn:C_R_gamma}
	C[R,\gamma] \to 0 \,\, \text{as} \,\, R\to 0.
\end{equation}

\vspace{.2in}
\noindent
Now, we choose the parameter $\gamma=\gamma_0$ such that the $C\big | \Psi ^i(0, \bm{0},\bm{0},\bm{0} ) \big | $ term in (\ref{eqn:theta_norm_bound_1}) satisfies
\begin{equation}
	\label{eqn:choose_gamma}
	C\big | \Psi ^i(0, \bm{0},\bm{0},\bm{0} ) \big | < \frac{\gamma _0}{4}.
\end{equation}

\vspace{.1in}
\noindent
Then, by (\ref{eqn:epsilon_R_gamma}) and (\ref{eqn:C_R_gamma}) we can choose $R$ sufficiently small such that the sum of the remaining terms in (\ref{eqn:theta_norm_bound_1}) is also less than $\displaystyle \frac{\gamma _0}{4} $. Therefore, we have proved that if  $ \displaystyle \|\bm{f}\|^{(2, \alpha)}  \leq \gamma _0,$ then for sufficiently small $R$,   $$ \|\Theta^i(\bm{f})\|^{(2, \alpha)}  \leq \frac{\gamma _0}{2}, $$ and this proves  $\bm{\Theta}$ maps $ \bm{\mathcal{E}}(R,\gamma _0)$ into itself.

\vspace{.2in}

\section {The Estimates for $ \| \bm{\Theta} (\bm{f}) - \bm{\Theta} (\bm{g}) \|^{(2, \alpha)}$}
\label{sec:estimate_f-g}

  \vspace{.2in}
  \noindent
  Now it remains to show that $  \bm{\Theta}$ is a contraction on $\bm{\mathcal{E}}(R, \gamma _0) $. The estimates are similar to those for $\displaystyle \| \bm{\Theta} (\bm{f})\|^{(2, \alpha)}$, so we will only point out the main steps without repeating all the calculations.
  
  \vspace{.2in}
  \noindent
  For any $\bm{f}, \bm{g} \in \bm{\mathcal{E}}(R, \gamma _0)$, from (\ref{eqn:theta}) we have  
  
  \allowdisplaybreaks
  \begin{equation}
  	\label{eqn:theta_norm_2}
  	\|\Theta^i(\bm{f}) - \Theta^i(\bm{g}) \|^{(2, \alpha)}\leq \|\omega^i(\bm{f}) - \omega^i(\bm{g})   \|^{(2,\alpha)}+\frac{1}{2}\Bigg | \sum_{\substack{k,l=1\\k\not =l }}^n\partial_k\partial_l\left (\omega^i(\bm{f}) - \omega^i(\bm{g} )  \right )(0) \Bigg |.
  \end{equation}
  
  \noindent
  \noindent
  By (\ref{defn:omega}) and Theorem \ref{thm: N(f)_bounded_by_f},
  
  \allowdisplaybreaks
  \begin{eqnarray}
  	\label{eqn: omega_estimate_2}
  	& & \|\omega^i(\bm{f})- \omega^i(\bm{g})  \|^{(2,\alpha)} \nonumber \\
  	& \leq & C \|\Psi^i \big (x, \bm{f}(x), D\bm{f}(x), D^2\bm{f}(x)\big ) - \Psi^i \big (x, \bm{g}(x), D\bm{g}(x), D^2\bm{g}(x)\big ) \|_{\alpha},  
  \end{eqnarray}
  
  \noindent
 and
 \allowdisplaybreaks
  \begin{eqnarray}
  	\label{eqn:omega_0_estimate_2}
  	& & \Big | \partial_k\partial_l\big (\omega^i(\bm{f}) - \omega^i(\bm{g})\big )(0) \Big | \nonumber \\
  	& \leq &   \sup_{B_R}\Big | \partial_k\partial_l\big (\omega^i(\bm{f}) - \omega^i(\bm{g})  \big ) \Big | \nonumber \\
  	& \leq & \|\omega^i(\bm{f}) - \omega^i(\bm{g}) \|^{(2,\alpha)} \nonumber \\
  	& \leq &  C \|\Psi^i \big (x,\bm{f}(x), D\bm{f}(x), D^2\bm{f}(x)\big ) - \Psi^i \big (x, \bm{g}(x), D\bm{g}(x), D^2\bm{g}(x)\big ) \|_{\alpha}. 
  \end{eqnarray}
  
  \vspace{.1in}
  \noindent
  Therefore we only need to estimate 
  
  $$\displaystyle  \|\Psi^i \big (x, \bm{f}(x), D\bm{f}(x), D^2\bm{f}(x)\big ) - \Psi^i \big (x, \bm{g}(x), D\bm{g}(x), D^2\bm{g}(x)\big )\|_{\alpha}.  $$
  
  \noindent
  Note that
  \allowdisplaybreaks
  \begin{eqnarray}
  	\label{eqn:psi_estimate_two}
  	&&\Psi^i \Big (x,\bm{f}(x), D\bm{f}(x), D^2\bm{f}(x) \Big ) - \Psi^i \Big (x,\bm{g}(x), D\bm{g}(x), D^2\bm{g}(x) \Big )     \nonumber\\
  	&=&\int_0^1 \frac{d}{dt}\Psi^i \Big ( x,t\bm{f}+(1-t)\bm{g},tD\bm{f}+(1-t)D\bm{g}, tD^2\bm{f}+(1-t)D^2\bm{g} \Big )dt\nonumber\\
  	&=& \sum_{j=1}^m A_j\left ( f^j - g^j \right )+\sum_{k=1}^n\sum_{j=1}^m B_k^j\partial_k\left ( f^j - g^j \right )+\sum_{j=1}^m \sum_{k,l=1}^nC_{kl}^j\partial_k\partial_l\left ( f^j - g^j \right ),
  \end{eqnarray}
  where $A_j$, $B_k^j$, and $C_{kl}^j$ are as defined in (\ref{eqn:const_deriv_psi_A}).

 \vspace{.2in}
 \noindent
 By Lemma \ref{lem:f-bounded-by_f(k)} and Lemma \ref{lem: compare_norms}, 
 
   \allowdisplaybreaks
 \begin{eqnarray}
 	\label{eqn:consant_domain_2}
	\addtolength{\jot}{2em}
 		\|f^j -g^j \|_{\alpha} & \leq & CR^2 \|f^j -g^j   \|^{(2, \alpha)}  \,\, \leq \,\, CR^2 \|\bm{f} -\bm{g} \|^{(2, \alpha)}, \nonumber \\
 		\|\partial _{k}(f^j -g^j  )\|_{\alpha} & \leq   & \|f^j-g^j \|^{(1, \alpha)} \,\, \leq \,\, CR \|f^j-g^j    \|^{(2, \alpha)} \,\, \leq \,\, CR \|\bm{f} -\bm{g} \|^{(2, \alpha)},  \\
 		\| \partial_{k}\partial_l (f^j -g^j )\|_{\alpha} & \leq & \|f^j -g^j \|^{(2, \alpha)} \,\, \leq \,\, \|\bm{f} -\bm{g} \|^{(2, \alpha)}. \nonumber
 \end{eqnarray}
 
 
 \vspace{.2in}
 \noindent
 Then by (\ref{eqn:psi_estimate_two}), (\ref{eqn:Banach_algebra}), (\ref{eqn:consant_domain_2}), (\ref{eqn:A_j_estimate}), (\ref{eqn:B_j_estimate}), and (\ref{eqn:C_j_estimate}) we have 
 
  \allowdisplaybreaks
 \begin{eqnarray}
 	& & \Big \|\Psi^i \Big (x,\bm{f}(x), D\bm{f}(x), D^2\bm{f}(x) \Big )  - \Psi^i \Big (x,\bm{g}(x), D\bm{g}(x), D^2\bm{g}(x)  \Big \|_{\alpha} \nonumber \\
 	& \leq & C \left ( R^2  \sum_{j=1}^m\|A_j\|_{\alpha} 
 	 +  R  \sum_{k=1}^n\sum_{j=1}^m \|B_k^j\|_{\alpha}  + \sum_{j=1}^m \sum_{k,l=1}^n\|C^j_{kl}\|_{\alpha}  \right ) \|\bm{f} -\bm{g} \|^{(2, \alpha)} \nonumber \\
 	 & \leq & C \left ( R^2  \sum_{j=1}^m\delta_A (R, \gamma _0) 
 	 +  R  \sum_{k=1}^n\sum_{j=1}^m \delta_B (R, \gamma _0)  + \sum_{j=1}^m \sum_{k,l=1}^n\delta_C (R, \gamma _0) \right ) \|\bm{f} -\bm{g} \|^{(2, \alpha)} \nonumber.	 
 \end{eqnarray}
 
 \vspace{.2in} 
  
 \noindent
 Note that by (\ref{eqn:delta_C_R_gamma}) and (\ref{eqn:C_R_gamma}),
 $$ 	\delta _C(R,\gamma _0) \to 0  \hspace{.2in} \text{as} \,\, R\to 0. $$
 
 
 \noindent
 Therefore  $$C \left ( R^2  \sum_{j=1}^m\delta_A (R, \gamma _0)  +  R  \sum_{k=1}^n\sum_{j=1}^m \delta_B (R, \gamma_0 )  + \sum_{j=1}^m \sum_{k,l=1}^n\delta_C (R, \gamma _0) \right ) \to 0 \hspace{.2in} \text{as} \,\, R \to 0. $$ 
 
 
 \vspace{.2in}
 \noindent
 Then by (\ref{eqn:theta_norm_2})-(\ref{eqn:omega_0_estimate_2}), when $R$ is sufficiently small, we have
 $$\|\Theta^i(\bm{f}) - \Theta^i(\bm{g}) \|^{(2, \alpha)}\leq \frac{1}{2} \|\bm{f} -\bm{g} \|^{(2, \alpha)},   $$ 
 
  \vspace{.1in}
 \noindent
 which implies $\bm{\Theta}$ is a contraction.   
 
 \vspace{.1in}
 \noindent
 Finally, by the Fixed Point Theorem there is a function $\bm{u}(x) \in \bm{\mathcal{E}}(R,\gamma _0)$ such that $$\bm{\Theta}(\bm{u})=\bm{u}. $$ Then 
 \begin{equation*}
 \Delta \bm{u} = - \bm{\Psi}\big (x, \bm{u}(x), D\bm{u}(x), D^2 \bm{u}(x) \big ),
 \end{equation*}

\vspace{.1in}
\noindent
with $\bm{u}(0)=\bm{0}$ and $D\bm{u}(0)=\bm{0}$.

\vspace{.1in}
\noindent
Thus Lemma \ref{thm: Poisson} is proved, and this completes the proof of Theorem \ref{thm: main}.

\vspace{.4in}

\bibliographystyle{plain}
\bibliography{thesis}

\end{document}